\documentclass[a4paper]{amsart}
\usepackage{graphicx}
\usepackage{amsmath}
\usepackage{amssymb}
\usepackage{oldgerm}
\usepackage[all]{xy}

\newcommand{\Hom}{\operatorname{Hom}\nolimits}

\renewcommand{\mod}{\operatorname{mod}\nolimits}

\newcommand{\id}{\operatorname{id}\nolimits}

\newcommand{\Ext}{\operatorname{Ext}\nolimits}

\newcommand{\HH}{\operatorname{HH}\nolimits}
\newcommand{\Id}{\operatorname{Id}\nolimits}
\newcommand{\Ho}{\operatorname{H}\nolimits}

\newcommand{\ra}{\mathfrak{r}}
\newcommand{\La}{\Lambda}

\newcommand{\cx}{\operatorname{cx}\nolimits}
\newcommand{\px}{\operatorname{px}\nolimits}
\newcommand{\op}{\operatorname{op}\nolimits}
\newcommand{\V}{\operatorname{V}\nolimits}

\newcommand{\T}{\operatorname{\mathcal{T}}\nolimits}
\newcommand{\C}{\operatorname{\mathcal{C}}\nolimits}
\newcommand{\D}{\operatorname{\mathcal{D}}\nolimits}

\newcommand{\e}{\operatorname{e}\nolimits}
\newcommand{\ev}{\operatorname{ev}\nolimits}
\newcommand{\Lae}{\Lambda^{\e}}

\newcommand{\MCM}{\operatorname{MCM}\nolimits}
\newcommand{\add}{\operatorname{add}\nolimits}
\newcommand{\thick}{\operatorname{thick}\nolimits}
\newcommand{\s}{\operatorname{\Sigma}\nolimits}
\newcommand{\Noeth}{\operatorname{Noeth}\nolimits}
\newcommand{\Noethfl}{\operatorname{Noeth^{\mathrm{fl}}}\nolimits}

\newtheorem{theorem}{Theorem}[section]

\newtheorem{corollary}[theorem]{Corollary}

\newtheorem{proposition}[theorem]{Proposition}
\theoremstyle{definition}
\newtheorem*{definition}{Definition}
\theoremstyle{definition}

\theoremstyle{definition}

\theoremstyle{definition}
\newtheorem*{example}{Example}
\theoremstyle{definition}

\theoremstyle{definition}

\theoremstyle{remark}

\theoremstyle{definition}

\theoremstyle{definition}

\begin{document}

\title{On the vanishing of cohomology in triangulated categories}
\author{Petter Andreas Bergh}
\address{Institutt for matematiske fag \\
  NTNU \\ N-7491 Trondheim \\ Norway}
\email{bergh@math.ntnu.no}


\thanks{The author was supported by NFR Storforsk grant no.\
167130}



\maketitle

\begin{abstract}
We study the vanishing of cohomology in a triangulated category, in
particular vanishing gaps and symmetry.
\end{abstract}

\section{Introduction}

In this paper, we study the vanishing of cohomology in a
triangulated category $\T$. Given two objects $X$ and $Y$ of $\T$, a
very natural question arises when looking at their cohomology: can
we detect the vanishing of $\Hom_{\T}(X, \s^nY)$ for large $n$ by
looking at \emph{finite} vanishing gaps? That is, is there a finite
set $S$ of integers such that the implication
$$\Hom_{\T}(X, \s^nY)=0 \text{ for } n \in S \hspace{2mm}
\Longrightarrow \hspace{2mm} \Hom_{\T}(X, \s^nY)=0 \text{ for } n
\gg 0$$ holds? Commutative local complete intersection rings provide
examples where this is true. Namely, for such a ring $A$, the
following was shown in \cite{Jorgensen} for a module $M$ of
complexity $d$: if $N$ is an $A$-module, and there exists an integer
$n> \dim A$ such that $\Ext_A^i(M,N)=0$ for $n \le i \le n+d$, then
$\Ext_A^i(M,N)$ vanishes for all $i> \dim A$. Over such a ring $A$,
the complexity of an $A$-module is at most the codimension $c$ of
$A$. Therefore, we can always detect vanishing of cohomology over
$A$ by looking at gaps of length $c+1$.

Another natural question is: does symmetry hold in the vanishing of
cohomology in $\T$? In other words, if $X$ and $Y$ are objects in
$\T$ such that $\Hom_{\T}(X, \s^nY)$ vanishes for $n \gg 0$, then
does it necessarily follow that $\Hom_{\T}(Y, \s^nX)$ also vanishes
for $n \gg 0$? Again, commutative local complete intersection rings
provide examples where this hold. Namely, it was shown in
\cite{AvramovBuchweitz} that if $M$ and $N$ are modules over such a
ring $A$, then the implication
$$\Ext_A^i(M,N)=0 \text{ for } i \gg 0 \hspace{2mm}
\Longrightarrow \hspace{2mm} \Ext_A^i(N,M)=0 \text{ for } i \gg 0$$
holds. Another class of rings where such symmetry holds are group
algebras of finite groups, or, more generally, as we shall see,
symmetric algebras with ``finitely generated" cohomology.

The two questions raised are studied in Section \ref{secvan} and
Section \ref{secsym}, respectively. We obtain some affirmative
answers when certain cohomology groups are finitely generated as a
module over a ring acting centrally on our triangulated category, a
concept we define in the following section. In particular, we show
that $\Ext$-symmetry holds for symmetric periodic algebras.

\section{Preliminaries}

Throughout this paper, we fix a triangulated category $\T$ with a
suspension functor $\s$. Thus $\T$ is an additive
$\mathbb{Z}$-category together with a class of distinguished
triangles satisfying Verdier's axioms (cf.\ \cite{Verdier}).

Recall that a \emph{thick} subcategory of $\T$ is a full
triangulated subcategory closed under direct summands. Now let $\C$
and $\D$ be subcategories of $\T$. We denote by $\thick^1_{\T} ( \C
)$ the full subcategory of $\T$ consisting of all the direct
summands of finite direct sums of shifts of objects in $\C$.
Furthermore, we denote by $\C \ast \D$ the full subcategory of $\T$
consisting of objects $M$ such that there exists a distinguished
triangle
$$C \to M \to D \to \s C$$
in $\T$, with $C \in \C$ and $D \in \D$. Now for each $n \ge 2$,
define inductively $\thick^n_{\T} ( \C )$ to be $\thick_{\T}^1 \left
( \thick^{n-1}_{\T} ( \C ) \ast \thick^1_{\T} ( \C ) \right )$, and
denote $\bigcup_{n=1}^{\infty} \thick^n_{\T} ( \C )$ by $\thick_{\T}
( \C )$. This is the smallest thick subcategory of $\T$ containing
$\C$.

The aim of this paper is to study the vanishing of cohomology in
triangulated categories satisfying a certain finite generation
hypothesis. This finite generation hypothesis is expressed in terms
of the \emph{graded center} $Z^*( \T )$ of our triangulated category
$\T$. Recall therefore that for an integer $n \in \mathbb{Z}$, the
degree $n$ component $Z^n( \T )$ is the set of natural
transformations $\Id \xrightarrow{f} \s^n$ satisfying $f_{\s X} =
(-1)^n \s f_X$ on the level of objects. For such a central element
$f$ and objects $X,Y \in \T$, consider the graded group
$\Hom_{\T}^*(X,Y) = \oplus_{i \in \mathbb{Z}} \Hom_{\T}(X, \s^i Y)$.
The element $f$ acts from the right on this graded group via the
morphism $X \xrightarrow{f_X} \s^n X$, and from the left via the
morphism $Y \xrightarrow{f_Y} \s^nY$. Namely, given a morphism $g
\in \Hom_{\T}(X, \s^m Y)$, the scalar product $gf$ is the
composition $X \xrightarrow{f_X} \s^n X \xrightarrow{\s^n g}
\s^{m+n}Y$, whereas $fg$ is the composition $X \xrightarrow{g} \s^m
Y \xrightarrow{\s^m f_Y} \s^{m+n} Y$. However, since $\Id
\xrightarrow{f} \s^n$ is a natural transformation, the diagram
$$\xymatrix{
X \ar[r]^g \ar[d]^{f_X} & \s^m Y \ar[d]^{f_{\s^m Y}} \\
\s^nX \ar[r]^<<<<<{\s^ng} & \s^{m+n}Y }$$ commutes, and so since
$f_{\s^mY}$ equals $(-1)^{mn} \s^mf_Y$ we see that $gf=(-1)^{mn}fg$.
This shows that $Z^*( \T )$ acts graded-commutatively on
$\Hom_{\T}^*(X,Y)$ for all objects $X$ and $Y$ in $\T$.

\sloppy Now let $R = \oplus_{i=0}^{\infty} R_i$ be a
graded-commutative ring, that is, for homogeneous elements $r_1,r_2
\in R$ the equality $r_1r_2 = (-1)^{|r_1||r_2|}r_2r_1$ holds. We say
that $R$ \emph{acts centrally} on $\T$ if there exists a graded ring
homomorphism $R \to Z^*( \T )$. If this is the case, then for every
object $X \in \T$ there is a graded ring homomorphism $R
\xrightarrow{\varphi_X} \Hom_{\T}^*(X,X)$ with the following
property: for all objects $Y \in \T$ the scalar actions from $R$ on
$\Hom_{\T}^*(X,Y)$ via $\varphi_X$ and $\varphi_Y$ are graded
equivalent, i.e.\
$$\varphi_Y(r) f = (-1)^{|r||f|} f \varphi_X(r)$$
for all homogeneous elements $r \in R$ and $f \in \Hom_{\T}^*(X,Y)$.
We say that the $R$-module $\Hom_{\T}^*(X,Y)$ is \emph{eventually
Noetherian}, and write $\Hom_{\T}^*(X,Y) \in \Noeth R$, if there
exists an integer $n_0 \in \mathbb{Z}$ such that the $R$-module
$\Hom_{\T}^{\ge n_0} (X,Y)$ is Noetherian. Moreover, we say that
$\Hom_{\T}^*(X,Y)$ is \emph{eventually Noetherian of finite length},
and write $\Hom_{\T}^*(X,Y) \in \Noethfl R$, if $\Hom_{\T}^*(X,Y)
\in \Noeth R$, and there exists an integer $n_0 \in \mathbb{Z}$ such
that $\ell_{R_0} \left ( \Hom_{\T} (X, \s^nY) \right ) < \infty$ for
each $n \ge n_0$. Note that if $\Hom_{\T}^*(X,Y)$ is eventually
Noetherian (respectively, eventually Noetherian of finite length),
then so is $\Hom_{\T}^*(X',Y')$ for all objects $X' \in
\thick_{\T}(X)$ and $Y' \in \thick_{\T}(Y)$. In particular, if our
category $\T$ is classically finitely generated, that is, if there
exists an object $G$ such that $\T = \thick_{\T}(G)$, then
$\Hom_{\T}^*(X,Y) \in \Noeth R$ (respectively, $\Hom_{\T}^*(X,Y) \in
\Noethfl R$) for all $X,Y \in \T$ if and only if $\Hom_{\T}^*(G,G)
\in \Noeth R$ (respectively, $\Hom_{\T}^*(G,G) \in \Noethfl R$).

\begin{definition}
Given objects $X$ and $Y$ of $\T$, we define the \emph{complexity}
of the ordered pair $(X,Y)$ as
$$\cx_{\T} (X,Y) \stackrel{\text{def}}{=} \dim_{R^{\ev}}
\Hom_{\T}^*(X,Y),$$ where $R^{\ev}$ denotes the commutative graded
subalgebra $\oplus_{i=0}^{\infty} R_{2i}$ of $R$. We define the
complexity $\cx_{\T}X$ of the single object $X$ as $\cx_{\T}X
\stackrel{\text{def}}{=} \cx_{\T} (X,X)$.
\end{definition}

When studying vanishing of cohomology in $\T$, we will only be
dealing with objects $X,Y \in \T$ with the property that the
$R$-module $\Hom_{\T}^*(X,Y)$ is eventually Noetherian of finite
length. This motivates the choice of terminology. Namely, it follows
from \cite[Proposition 2.6]{BIKO} that if $\Hom_{\T}^*(X,Y) \in
\Noethfl R$, then the Krull dimension  of the $R^{\ev}$-module
$\Hom_{\T}^*(X,Y)$ equals the infimum of all non-negative integers
$t$ with the following property: there exists a real number $a$ such
that
$$\ell_{R_0} \left ( \Hom_{\T}(X, \s^nY) \right ) \le an^{t-1}$$
for $n \gg 0$. A priori, the complexity of a pair is not finite.
However, when $\Hom_{\T}^*(X,Y)$ is eventually Noetherian of finite
length, then the finiteness of $\cx_{\T} (X,Y)$ follows from the
above together with \cite[Remark 2.1]{BIKO} and \cite[Theorem
11.1]{AtiyahMacdonald}.

It follows from the above alternative description of complexity that
if $X$ and $Y$ are objects of $\T$ with $\Hom_{\T}^*(X,Y) \in
\Noethfl R$, then $\cx_{\T}(X,Y) =0$ if and only if
$\Hom_{\T}^*(X,Y)$ is eventually zero, that is, if $\Hom_{\T}(X,
\s^nY) =0$ for $n \gg 0$. Now digress for a moment, and let $\La$ be
a ring. Then $\La$ satisfies \emph{Auslander's condition} if for
every finitely generated module $M$, there exists an integer $d_M$,
depending only on $M$, satisfying the following: if $N$ is a
finitely generated $\La$-module and $\Ext_{\La}^n(M,N)=0$ for $n \gg
0$, then $\Ext_{\La}^n(M,N)=0$ for $n \ge d_M$. Motivated by this,
we define a full subcategory $\C$ of $\T$ to be a \emph{left
Auslander subcategory} if for every object $X \in \C$, there exists
an integer $d_X$, depending only on $X$, such that the following
holds: if $\Hom_{\T}^*(X,Y)$ is eventually zero for some object $Y
\in \T$, then $\Hom_{\T}(X, \s^nY) =0$ for $n \ge d_X$. It is easy
to see that this holds if and only if for all objects $X \in \C$ and
$Y \in \T$, the implication
$$\Hom_{\T}(X, \s^n Y) =0 \text{ for } n \gg 0 \hspace{2mm}
\Longrightarrow \hspace{2mm} \Hom_{\T}(X, \s^n Y)=0 \text{ for all }
n \in \mathbb{Z}$$ holds. Dually, we can define right Auslander
subcategories. Note that if an object $X \in \T$ belongs to a left
or right Auslander subcategory and $\Hom_{\T}^*(X,X)$ is eventually
zero, then $X=0$.

\section{Vanishing of cohomology}\label{secvan}

We start with the following result, the key ingredient in the main
theorem. It shows that we can always reduce the complexity of an
object whose endomorphism ring is eventually Noetherian of finite
length. However, recall first the following notion. Let $R$ be a
graded-commutative ring acting centrally on $\T$, and let $X \in \T$
be an object. Then, given a homogeneous element $r \in R$, we can
complete the map $X \xrightarrow{\varphi_X(r)} \s^{|r|} X$ into a
triangle
$$X \xrightarrow{\varphi_X(r)} \s^{|r|} X \to X/\!\!/r \to \s X.$$
The object $X /\!\!/r$ is well defined up to isomorphism, and is
called a \emph{Koszul object} of $r$ on $X$.

\begin{proposition}\label{redcx}
Let $R$ be a graded-commutative ring acting centrally on $\T$, and
let $X \in \T$ be an object such that $\Hom_{\T}^*(X,X) \in \Noethfl
R$. Then if $\cx_{\T}X$ is nonzero, there exists a homogeneous
element $r \in R$, of positive degree, whose Koszul object
$X/\!\!/r$ in the triangle
$$X \xrightarrow{\varphi_X(r)} \s^{|r|} X \to X/\!\!/r \to \s X$$
satisfies $\cx_{\T}X/\!\!/r = \cx_{\T}X-1$.
\end{proposition}

\begin{proof}
Suppose $\cx_{\T}X >0$. By \cite[Lemma 2.5]{BIKO}, there exists an
integer $i_0$ and a homogeneous element $r \in R$, of positive
degree, such that scalar multiplication
$$\Hom_{\T}(X,\s^i X) \xrightarrow{\cdot r} \Hom_{\T}(X,\s^{i+|r|} X)$$
is injective for $i \ge i_0$. Applying $\Hom_{\T}(X,-)$ to the triangle
$$X \xrightarrow{\varphi_X(r)} \s^{|r|} X \to X/\!\!/r \to \s X,$$
we obtain a long exact sequence
$$\cdots \to \Hom_{\T}(X, \s^i X) \xrightarrow{\cdot (-1)^i r}
\Hom_{\T}(X,\s^{i+|r|} X) \to \Hom_{\T}(X, \s^i X/\!\!/r ) \to
\cdots$$ in cohomology. This long exact sequence induces a short
exact sequence
$$0 \to \Hom_{\T}^{\ge i_0}(X,X) \xrightarrow{\cdot r} \Hom_{\T}^{\ge
i_0 + |r|}(X,X) \to \Hom_{\T}^{\ge i_0}(X,X/\!\!/r) \to 0$$ of
eventually Noetherian $R$-modules of finite length, a sequence from
which we deduce that $\cx_{\T}(X,X/\!\!/r) = \cx_{\T}X-1$. From this
exact sequence we also see that the element $r$ annihilates
$\Hom_{\T}^{\ge i_0}(X,X/\!\!/r)$. Therefore, when applying
$\Hom_{\T}(-,X/\!\!/r)$ to our triangle, we obtain a short exact
sequence
$$0 \to  \Hom_{\T}^{\ge i_0+|r|}(X,X/\!\!/r) \to \Hom_{\T}^{\ge
i_0+|r|+1}(X/\!\!/r,X/\!\!/r) \to \Hom_{\T}^{\ge i_0+1}(X,X/\!\!/r)
\to 0$$ of eventually Noetherian $R$-modules of finite length. This
gives the inequality $\cx_{\T} X/\!\!/r \le \cx_{\T}(X,X/\!\!/r)$.

Since $\Hom_{\T}^*(X,X/\!\!/r)$ is eventually Noetherian, there
exists an integer $n_0$ such that the $R$-module $\Hom_{\T}^{\ge
n_0}(X,X/\!\!/r)$ is finitely generated. The $R$-scalar action
factors through the ring homomorphism $\varphi_{X/\!\!/r}$, and
therefore $\Hom_{\T}^{\ge n_0}(X,X/\!\!/r)$ is also finitely
generated as a module over $\Hom_{\T}^*(X/\!\!/r,X/\!\!/r)$.
Consequently, the rate of growth of the sequence $\{ \ell_{R_0}
\left ( \Hom_{\T}(X, \s^n X/\!\!/r) \right ) \}_{n=1}^{\infty}$ is
at most the rate of growth of $\{ \ell_{R_0} \left (
\Hom_{\T}(X/\!\!/r, \s^n X/\!\!/r) \right ) \}_{n=1}^{\infty}$,
i.e.\ $\cx_{\T}(X,X/\!\!/r) \le \cx_{\T} X/\!\!/r$.
\end{proof}

Now let $X$ be an object of $\T$, and let $r_1, \dots, r_c$ be a
sequence of homogeneous elements belonging to a graded-commutative
ring $R$ acting centrally on $\T$, with $c \ge 2$. Then we define
the Koszul element $X/\!\!/(r_1, \dots, r_c)$ inductively as
$$X/\!\!/(r_1, \dots, r_c) \stackrel{\text{def}}{=} \left (
X/\!\!/(r_1, \dots, r_{c-1}) \right ) /\!\!/r_c.$$ Suppose the
$R$-module $\Hom_{\T}^*(X,X)$ is eventually Noetherian of finite
length, and denote the complexity of $X$ by $c$. If $c>0$, then the
previous result guarantees the existence of a homogeneous element
$r_1 \in R$, of positive degree, such that $\cx_{\T} X/\!\!/r_1
=c-1$. Since $X/\!\!/r_1$ belongs to $\thick_{\T}(X)$, the
$R$-module $\Hom_{\T}^*(X/\!\!/r_1,X/\!\!/r_1)$ is also eventually
Noetherian of finite length. Therefore, if $c-1>0$, then we may use
the above result again; there exists a homogeneous element $r_2 \in
R$, of positive degree, such that $\cx_{\T} X/\!\!/(r_1,r_2) =c-2$.
Continuing like this, we obtain a sequence $r_1, \dots, r_c$ of
homogeneous elements of $R$, all of positive degree, together with
triangles
\begin{eqnarray*}
& X \to \s^{|r_1|} X \to X_1 \to \s X & \\
& X_1 \to \s^{|r_2|} X_1 \to X_2 \to \s X_1 & \\
& \vdots & \\
& X_{c-1} \to \s^{|r_c|} X_{c-1} \to X_c \to \s X_{c-1} &
\end{eqnarray*}
in which $X_i = X/\!\!/(r_1, \dots, r_i)$ and $\cx_{\T} X_i = c-i$
for $1 \le i \le c$. We say that the sequence $r_1, \dots, r_c$
\emph{reduces the complexity} of the object $X$. Note that such a
sequence is not unique in general, and that the order matters. Note
also that the sequence reduces the complexity of $\s^n X$ for any $n
\in \mathbb{Z}$, since $\cx_{\T}Y = \cx_{\T} \s^nY$ for all objects
$Y$ in $\T$.

We now prove our main result. It shows that, for two objects $X$ and
$Y$, if $\Hom_{\T}^*(X,Y)$ contains a large enough ``gap", then
$\Hom_{\T}^*(X,Y)$ is actually zero. The length of the gap depends
on the sum of the degrees of a sequence reducing the complexity of
$X$.

\begin{theorem}\label{mainthm1}
Let $X$ and $Y$ be objects of $\T$ with $\Hom_{\T}^*(X,X) \in
\Noethfl R$ for some graded-commutative ring $R$ acting centrally on
$\T$, and suppose $\thick_{\T} (X)$ is either a left or right
Auslander subcategory of $\T$. Let $r_1, \dots, r_c$ be a sequence
of positive degree homogeneous elements of $R$ reducing the
complexity of $X$, where $c = \cx_{\T}X$. Then the following are
equivalent:
\begin{enumerate}
\item[(i)] There exists an integer $n \in \mathbb{Z}$ such that
$\Hom_{\T}(X, \s^i Y) =0$ for $n \le i \le n+|r_1|+ \cdots
+|r_c|-c$.
\item[(ii)] $\Hom_{\T}(X, \s^i Y) =0$ for all $i \in \mathbb{Z}$.
\end{enumerate}
\end{theorem}

\begin{proof}
We argue by induction on $c$ that (i) implies (ii). If $c$ is zero,
then $\Hom_{\T}^*(X,X)$ is eventually zero, and so $X=0$ since
$\thick_{\T} (X)$ is either a left or right Auslander subcategory of
$\T$. If $c>0$, consider the triangle
$$X \to \s^{|r_1|}X \to X/\!\!/r_1 \to \s X.$$
Applying $\Hom_{\T}(-,Y)$ to this triangle gives the long exact sequence
$$\cdots \to \Hom_{\T}(X, \s^{i-1}Y) \to \Hom_{\T}(X/\!\!/r_1, \s^iY)
\to \Hom_{\T}(X, \s^{i-|r_1|}Y) \to \cdots$$ in cohomology, from
which we see that $\Hom_{\T}(X/\!\!/r_1, \s^iY)=0$ for
$$(n+|r_1|) \le i \le (n+|r_1|)+ |r_2|+ \cdots +|r_c|-(c-1).$$
The complexity of $X/\!\!/r_1$ is $c-1$, and the sequence $r_2,
\dots, r_c$ is a complexity reducing sequence for this object.
Therefore, by the induction hypothesis, we conclude that
$\Hom_{\T}(X/\!\!/r_1, \s^iY)=0$ for all $i \in \mathbb{Z}$. The
long exact sequence then shows that $\Hom_{\T}(X, \s^i Y)$ is
isomorphic to $\Hom_{\T}(X, \s^{i+|r_1|}Y)$ for all integers $i$,
and from (i) we then see that $\Hom_{\T}(X, \s^i Y)$ must vanish for
all $i \in \mathbb{Z}$.
\end{proof}

If we interchange the objects $X$ and $Y$ in the theorem, then the
corresponding result of course holds. We state this without proof.

\begin{theorem}\label{mainthm2}
Let $X$ and $Y$ be objects of $\T$ with $\Hom_{\T}^*(X,X) \in
\Noethfl R$ for some graded-commutative ring $R$ acting centrally on
$\T$, and suppose $\thick_{\T} (X)$ is either a left or right
Auslander subcategory of $\T$. Let $r_1, \dots, r_c$ be a sequence
of positive degree homogeneous elements of $R$ reducing the
complexity of $X$, where $c = \cx_{\T}X$. Then the following are
equivalent:
\begin{enumerate}
\item[(i)] There exists an integer $n \in \mathbb{Z}$ such that
$\Hom_{\T}(Y, \s^i X) =0$ for $n \le i \le n+|r_1|+ \cdots
+|r_c|-c$.
\item[(ii)] $\Hom_{\T}(Y, \s^i X) =0$ for all $i \in \mathbb{Z}$.
\end{enumerate}
\end{theorem}

We now use these results to study the vanishing of cohomology over
Artin algebras. Let $\La$ be a Noetherian ring, and denote the
bounded derived category of finitely generated $\La$-modules by
$D^b( \La )$. Furthermore, let $D^{\text{perf}} ( \La )$ be the
thick subcategory of $D^b( \La )$ generated by $\La$; it consists of
the perfect complexes, that is, objects isomorphic to bounded
complexes of finitely generated projective $\La$-modules. The
\emph{stable derived category} of $\La$, denoted $D^b_{\text{st}} (
\La )$, is the Verdier quotient
$$D^b_{\text{st}} ( \La ) \stackrel{\text{def}}{=}  D^b( \La )/
D^{\text{perf}} ( \La ).$$ This is a triangulated category whose
suspension functor corresponds to that in $D^b( \La )$. Moreover,
by \cite[Remark 5.1]{BIKO}, the central action of a
graded-commutative ring $R$ on $D^b( \La )$ carries over to
$D^b_{\text{st}} ( \La )$ via the ring homomorphism $Z^*(D^b( \La
)) \to Z^*(D^b_{\text{st}} ( \La ))$ induced by the natural
quotient functor. Thus if $X$ and $Y$ are complexes in $D^b( \La
)$, then the natural map
$$\Hom_{D^b( \La )}^*(X,Y) \to \Hom_{D^b_{\text{st}}( \La )}^*(X,Y)$$
is an $R$-module homomorphism. If $\La$ is also Gorenstein, that is,
if the injective dimension of $\La$ both as a left and as a right
module over itself is finite, then by \cite[Corollary
6.3.4]{Buchweitz} this homomorphism is eventually bijective. That
is, if $\La$ is a Noetherian Gorenstein ring, then the natural map
$$\Hom_{D^b( \La )}(X, \s^n Y) \to \Hom_{D^b_{\text{st}}( \La )}(X, \s^n Y)$$
is bijective for $n \gg 0$. Consequently, for any complexes $X$ and
$Y$ in $D^b( \La )$, if $\Hom_{D^b( \La )}^*(X,Y)$ is an eventually
Noetherian $R$-module of finite length, then so is
$\Hom_{D^b_{\text{st}}( \La )}^*(X,Y)$.

Suppose $\La$ is an Artin algebra, that is, the center $Z( \La )$ of
$\La$ is a commutative Artin ring over which $\La$ is finitely
generated as a module. Denote by $\mod \La$ the category of finitely
generated left $\La$-modules. If $\La$ is Gorenstein, then denote by
$\MCM ( \La )$ the category of finitely generated maximal
Cohen-Macaulay $\La$-modules, i.e.\
$$\MCM ( \La ) = \{ M \in \mod \La \mid \Ext_{\La}^i(M, \La )=0
\text{ for all } i >0 \}.$$ It follows from general cotilting
theory that this is a Frobenius exact category, in which the
projective injective objects are the projective $\La$-modules, and
the injective envelopes are the left $\add \La$-approximations.
Therefore the stable category $\underline{\MCM} ( \La )$, which is
obtained by factoring out all morphisms which factor through
projective $\La$-modules, is a triangulated category. Its shift
functor is given by cokernels of left $\add \La$-approximations,
the inverse shift is the usual syzygy functor. It follows from
work by Buchweitz, Happel and Rickard (cf.\ \cite{Buchweitz},
\cite{Happel}, \cite{Rickard}) that $\underline{\MCM} ( \La )$ and
the quotient category $D^b ( \La ) / D^{\text{perf}} ( \La )$ are
equivalent as triangulated categories. If $M$ and $N$ are maximal
Cohen-Macaulay modules in $\mod \La$, then there is an isomorphism
$$\Ext_{\La}^n(M,N) \simeq \Hom_{\underline{\MCM}( \La )}(
\Omega_{\La}^n(M),N)$$ for every $n>0$. We use this isomorphism to
prove the following result. It shows that when a certain finiteness
condition holds, then the thick subcategory in $\underline{\MCM}(
\La )$ generated by a module is a left and right Auslander
subcategory.

\begin{proposition}\label{auscat}
Let $\La$ be an Artin Gorenstein algebra with Jacobson radical
$\ra$, and let $M$ be a maximal Cohen-Macaulay module. If either
$\Ext_{\La}^*(M, \La / \ra )$ or $\Ext_{\La}^*( \La / \ra, M )$
belongs to $\Noethfl R$ for some graded-commutative ring $R$ acting
centrally on $D^b ( \La )$, then $\thick_{\underline{\MCM}( \La
)}(M)$ is a left and right Auslander subcategory of
$\underline{\MCM}( \La )$.
\end{proposition}

\begin{proof}
Suppose that $\Ext_{\La}^*(M, \La / \ra ) \in \Noethfl R$. Let $X$
and $Y$ be maximal Cohen-Macaulay modules in $\mod \La$ with $X \in
\thick_{\underline{\MCM}( \La )}(M)$, and suppose that
$\Hom_{\underline{\MCM}( \La )}^*(X,Y)$ is eventually zero. We prove
by induction on $\cx_{\underline{\MCM}( \La )}X$ that
$\Hom_{\underline{\MCM}( \La )}(X, \s^nY)=0$ for all $n \in
\mathbb{Z}$.

Suppose $\cx_{\underline{\MCM}( \La )}X=0$. The finiteness condition
implies that the $R$-module $\Ext_{\La}^*(M,N)$ is eventually
Noetherian of finite length for all $N \in \mod \La$, and so the
same holds for $\Hom_{\underline{\MCM}( \La )}^*(X,X)$. Therefore
$\Hom_{\underline{\MCM}( \La )}^*(X,X)$ is eventually zero, in
particular $\Ext_{\La}^n(X,X) =0$ for $n \gg 0$. Now consider the
$R$-module $\Ext_{\La}^*(X, \La / \ra )$. Since it belongs to
$\Noethfl R$, and the $R$-module structure factors through the ring
homomorphism $R \xrightarrow{\varphi_X} \Ext_{\La}^*(X,X)$, we see
that it must be eventually zero. The $\La$-module $X$ therefore has
finite projective dimension, and is isomorphic to the zero object in
$\underline{\MCM}( \La )$. Consequently $\Hom_{\underline{\MCM}( \La
)}(X, \s^nY)=0$ for all $n \in \mathbb{Z}$.

If $\cx_{\underline{\MCM}( \La )}X>0$, then let $r \in R$ be a
homogeneous element of positive degree such that
$\cx_{\underline{\MCM}( \La )}X/\!\!/r = \cx_{\underline{\MCM}( \La
)}X -1$. Since $\Hom_{\underline{\MCM}( \La )}^*(X,Y)$ is eventually
zero, we see from the triangle
$$X \to \s^{|r|} X \to X/\!\!/r \to \s X$$
that the same holds for $\Hom_{\underline{\MCM}( \La
)}^*(X/\!\!/r,Y)$. The Koszul object $X/\!\!/r$ belongs to
$\thick_{\underline{\MCM}( \La )}(M)$, hence by induction
$\Hom_{\underline{\MCM}( \La )}(X/\!\!/r, \s^nY)=0$ for all $n \in
\mathbb{Z}$. From the triangle we obtain the isomorphism
$$\Hom_{\underline{\MCM}( \La )}(X, \s^nY) \simeq \Hom_{\underline{\MCM}( \La )}(X,
\s^{n+|r|}Y)$$ for all integers $n$, and this implies that
$\Hom_{\underline{\MCM}( \La )}(X, \s^nY)=0$ for all $n \in
\mathbb{Z}$.

We have now proved that if $\Ext_{\La}^*(M, \La / \ra ) \in \Noethfl
R$, then $\thick_{\underline{\MCM}( \La )}(M)$ is a left Auslander
subcategory of $\underline{\MCM}( \La )$. Virtually the same proof
shows that $\thick_{\underline{\MCM}( \La )}(M)$ is also a right
Auslander subcategory of $\underline{\MCM}( \La )$. Moreover, an
analogous proof shows that the same holds if $\Ext_{\La}^*( \La /
\ra, M ) \in \Noethfl R$.
\end{proof}

Before proving the next result, we recall the following. Let $\La$
be an Artin algebra with Jacobson radical $\ra$, and let $M \in \mod
\La$ be a module with minimal projective and injective resolutions
$$\cdots \to P_2 \to P_1 \to P_0 \to M \to 0$$
and
$$0 \to M \to I^0 \to I^1 \to I^2 \to \cdots$$
respectively. Then the \emph{complexity} and \emph{plexity} of $M$,
denoted $\cx_{\La} M$ and $\px_{\La} M$, respectively, are defined
as
\begin{eqnarray*}
\cx_{\La} M & \stackrel{\text{def}}{=} & \inf \{ t \in \mathbb{N}
\cup \{ 0 \} \mid \exists a \in \mathbb{R} \text{ such that }
\ell_{Z(
\La )} (P_n) \le an^{t-1} \text{ for } n \gg 0 \}, \\
\px_{\La} M & \stackrel{\text{def}}{=} & \inf \{ t \in \mathbb{N}
\cup \{ 0 \} \mid \exists a \in \mathbb{R} \text{ such that }
\ell_{Z( \La )} (I^n) \le an^{t-1} \text{ for } n \gg 0 \},
\end{eqnarray*}
where $Z( \La )$ is the center of $\La$. Now let $R$ be a
graded-commutative ring acting centrally on $D^b ( \La )$. If
$\Ext_{\La}^*(M, \La / \ra )$ belongs to $\Noethfl R$, then
$\cx_{\La} M$ coincides with $\cx_{D^b( \La )}M$, and
$\Ext_{\La}^*(M, M)$ also belongs to $\Noethfl R$. Therefore there
exists a sequence $r_1, \dots, r_{\cx_{\La} M}$ of homogeneous
elements of $R$, all of positive degree, reducing the complexity of
$M$ as an object in $D^b( \La )$. Similarly, if $\Ext_{\La}^*( \La /
\ra, M )$ belongs to $\Noethfl R$, then $\px_{\La} M = \cx_{D^b( \La
)}M$. In this case, there exists a homogeneous sequence in $R$ of
length $\px_{\La} M$ reducing the complexity of $M$ as an object in
$D^b( \La )$. Using this and Proposition \ref{auscat}, we obtain the
following vanishing results on cohomology over Gorenstein algebras.
We prove only the first of these results; the proof of the other
result is similar.

\begin{theorem}\label{vangor1}
Let $\La$ be an Artin Gorenstein algebra with Jacobson radical
$\ra$, and let $M \in \mod \La$ be a maximal Cohen-Macaulay module.
Suppose $\Ext_{\La}^*(M, \La / \ra ) \in \Noethfl R$ for some
graded-commutative ring $R$ acting centrally on $D^b ( \La )$, and
let $r_1, \dots, r_c$ be a sequence of positive degree homogeneous
elements of $R$ reducing the complexity of $M$, where $c = \cx_{\La}
M$. Then for any $N \in \mod \La$, the implications \emph{(i)}
$\Leftrightarrow$ \emph{(ii)} and \emph{(iii)} $\Leftrightarrow$
\emph{(iv)} hold for the following statements:
\begin{enumerate}
\item[(i)] There exists a number $n > \id \La$ such that
$\Ext_{\La}^i(M,N)=0$ for $n \le i \le n+|r_1|+ \cdots +|r_c|-c$.
\item[(ii)] $\Ext_{\La}^i(M,N)=0$ for all $i > \id \La$.
\item[(iii)] There exists a number $n > \id \La$ such that
$\Ext_{\La}^i(N,M)=0$ for $n \le i \le n+|r_1|+ \cdots +|r_c|-c$.
\item[(iv)] $\Ext_{\La}^i(N,M)=0$ for all $i > \id \La$.
\end{enumerate}
\end{theorem}

\begin{proof}
By \cite[Theorem 1.8]{AuslanderBuchweitz}, there exists an exact
sequence
$$0 \to Q \to C \to N \to 0$$
in $\mod \La$, in which $Q$ has finite projective dimension and $C$
is maximal Cohen-Macaulay. Since $Q$ also has finite injective
dimension and $\id Q$ is at most $\id \La$, there are isomorphisms
$\Ext_{\La}^i(M,N) \simeq \Ext_{\La}^i(M,C)$ for $i > \id \La$.
Moreover, the module $\Omega_{\La}^{\id \La}(N)$ is maximal
Cohen-Macaulay, and $\Ext_{\La}^i(N,M) \simeq \Ext_{\La}^{i- \id
\La}( \Omega_{\La}^{\id \La}(N),M)$ for $i > \id \La$. We may
therefore without loss of generality assume that $N$ itself is
maximal Cohen-Macaulay, and replace $\id \La$ by $0$ in the
statements. The implications now follow from Theorem \ref{mainthm1},
Theorem \ref{mainthm2}, Proposition \ref{auscat} and the fact that
$\Ext_{\La}^i(X,Y) \simeq \Hom_{\underline{\MCM}( \La )}(X, \s^iY)$
when $X$ and $Y$ are maximal Cohen-Macaulay and $i>0$.
\end{proof}

\begin{theorem}\label{vangor2}
Let $\La$ be an Artin Gorenstein algebra with Jacobson radical
$\ra$, and let $M \in \mod \La$ be a maximal Cohen-Macaulay module.
Suppose $\Ext_{\La}^*( \La / \ra, M ) \in \Noethfl R$ for some
graded-commutative ring $R$ acting centrally on $D^b ( \La )$, and
let $r_1, \dots, r_c$ be a sequence of positive degree homogeneous
elements of $R$ reducing the complexity of $M$, where $c = \px_{\La}
M$. Then for any $N \in \mod \La$, the implications \emph{(i)}
$\Leftrightarrow$ \emph{(ii)} and \emph{(iii)} $\Leftrightarrow$
\emph{(iv)} hold for the following statements:
\begin{enumerate}
\item[(i)] There exists a number $n > \id \La$ such that
$\Ext_{\La}^i(M,N)=0$ for $n \le i \le n+|r_1|+ \cdots +|r_c|-c$.
\item[(ii)] $\Ext_{\La}^i(M,N)=0$ for all $i > \id \La$.
\item[(iii)] There exists a number $n > \id \La$ such that
$\Ext_{\La}^i(N,M)=0$ for $n \le i \le n+|r_1|+ \cdots +|r_c|-c$.
\item[(iv)] $\Ext_{\La}^i(N,M)=0$ for all $i > \id \La$.
\end{enumerate}
\end{theorem}

\section{Symmetry}\label{secsym}

For a group algebra $kG$ of a finite group $G$ over a field $k$, the
group cohomology ring $\Ho^*(G,k)$ is graded-commutative and acts
centrally on $D^b( kG)$. Moreover, by a classical result of Evens
and Venkov (cf.\ \cite{Evens}, \cite{Venkov1}, \cite{Venkov2}), the
cohomology ring is Noetherian, and $\Ext_{kG}^*(M,N)$ is a finitely
generated $\Ho^*(G,k)$-module for all $M$ and $N$ in $\mod kG$.
Therefore the vanishing results from the previous section apply to
group algebras.

Commutative local complete intersection rings also have finitely
generated cohomology. For such a ring $A$, it was shown in
\cite{Avramov1} that there exists a certain polynomial ring
$\widehat{A}[\chi_1, \dots, \chi_c]$ acting centrally on $D^b (
\widehat{A} )$, where $\widehat{A}$ denotes the completion of $A$
with respect to its maximal ideal. Again, for all finitely generated
$A$-modules $M$ and $N$, the $\widehat{A}[\chi_1, \dots,
\chi_c]$-module $\Ext_{\widehat{A}}^*( \widehat{M}, \widehat{N} )$
is finitely generated. Consequently, vanishing results similar to
those in the previous section also hold in this case.

A fascinating aspect of the vanishing of cohomology over both group
algebras and commutative local complete intersections is symmetry.
In \cite{AvramovBuchweitz} it was shown that for finitely generated
modules $M$ and $N$ over a commutative local complete intersection
$A$, the vanishing of $\Ext_A^i (M,N)$ for $i \gg 0$ implies the
vanishing of $\Ext_A^i (N,M)$ for $i \gg 0$. The proof involves the
theory of certain support varieties attached to each pair of
$A$-modules. Denote by $c$ the codimension of $A$ and by $K$ the
algebraic closure of its residue field. A cone $\V^*_A (M,N)$ in
$K^c$ is associated to the ordered pair $(M,N)$, with the following
properties:
\begin{eqnarray*}
\V^*_A (M,N)= \{ 0 \} & \Leftrightarrow & \Ext_A^i (M,N)=0 \text{
for } i \gg 0, \\
\V^*_A (M,N) & = & \V^*_A (M,M) \cap \V^*_A (N,N).
\end{eqnarray*}
The symmetry in the vanishing of cohomology follows immediately from
these properties. Similarly, the theory of support varieties for
modules over group algebras of finite groups can be used to show
that symmetry holds also for such algebras (cf.\ \cite{Benson}).

We shall see in this section that in general there is no symmetry in
the vanishing of cohomology over an Artin algebra, even when the
algebra is selfinjective and has finitely generated cohomology in
the sense of group algebras. But first, we study situations where
symmetry holds. Let $k$ be a commutative Artin ring, and suppose
$\T$ is a $\Hom$-finite triangulated $k$-category. In other words,
for all objects $X,Y,Z \in \T$ the group $\Hom_{\T}(X,Y)$ is a
$k$-module of finite length, and composition
$$\Hom_{\T}(Y,Z) \times \Hom_{\T}(X,Y) \to \Hom_{\T}(X,Z)$$
is $k$-bilinear, where $D = \Hom_k(-,k)$. A \emph{Serre functor} on
$\T$ is a triangle equivalence $\T \xrightarrow{S} \T$, together
with functorial isomorphisms
$$\Hom_{\T}(X,Y) \simeq D \Hom_{\T}(Y,SX)$$
of $k$-modules for all objects $X,Y \in \T$. By
\cite{BondalKapranov}, such a functor is unique if it exists.
Following \cite{Keller}, for an integer $d \in \mathbb{Z}$, the
category $\T$ is said to be \emph{weakly $d$-Calabi-Yau} if it
admits a Serre functor which is isomorphic as a $k$-linear functor
to $\s^d$. If, in addition, this isomorphism is an isomorphism of
triangle functors, then $\T$ is \emph{$d$-Calabi-Yau}. However, we
will only be dealing with weakly $d$-Calabi-Yau categories. When
$\T$ is such a category, then for all objects $X,Y \in \T$ there is
an isomorphism
$$\Hom_{\T}(X,Y) \simeq D \Hom_{\T}(Y, \s^dY)$$
of $k$-modules. It follows immediately that if this holds, then
$\Hom_{\T}(X, \s^nY)=0$ for $n \gg 0$ if and only if $\Hom_{\T}(Y,
\s^nX)=0$ for $n \ll 0$.

Now let $\La$ be an Artin Gorenstein algebra. Following \cite{Mori},
we say that $\La$ is \emph{stably symmetric} if $\underline{\MCM}(
\La )$ is weakly $d$-Calabi-Yau for some integer $d \in \mathbb{Z}$.
It was shown in that paper that if $\La$ in addition satisfies
Auslander's condition, then symmetry holds in the vanishing of
cohomology of $\La$-modules. The following result shows that
symmetry holds for modules with finitely generated cohomology.

\begin{theorem}\label{localsym}
Let $\La$ be a stably symmetric Artin Gorenstein algebra with
Jacobson radical $\ra$. Let $M \in \mod \La$ be a module such that
either $\Ext_{\La}^*(M, \La / \ra )$ or $\Ext_{\La}^*( \La / \ra, M
)$ belongs to $\Noethfl R$ for some graded-commutative ring $R$
acting centrally on $D^b ( \La )$. Then for every $N \in \mod \La$,
the following are equivalent:
\begin{enumerate}
\item[(i)] $\Ext_{\La}^i(M,N)=0$ for $i \gg 0$.
\item[(ii)] $\Ext_{\La}^i(M,N)=0$ for $i > \id \La$.
\item[(iii)] $\Ext_{\La}^i(N,M)=0$ for $i \gg 0$.
\item[(iv)] $\Ext_{\La}^i(N,M)=0$ for $i > \id \La$.
\end{enumerate}
\end{theorem}

\begin{proof}
As in the proof of Theorem \ref{vangor1}, there exists an exact
sequence
$$0 \to Q_N \to C_N \to N \to 0$$
in $\mod \La$, in which $Q_N$ has finite projective (and injective)
dimension, and $C_N$ is maximal Cohen-Macaulay. Thus there is an
isomorphism $\Ext_{\La}^i(M,N) \simeq \Ext_{\La}^{i- \id \La}(
\Omega_{\La}^{\id \La}(M),C_N)$ for every $i> \id \La$. Moreover,
since either $\Ext_{\La}^*(M, \La / \ra )$ or $\Ext_{\La}^*( \La /
\ra, M )$ belongs to $\Noethfl R$, so do either $\Ext_{\La}^*(
\Omega_{\La}^{\id \La}(M), \La / \ra )$ or $\Ext_{\La}^*( \La / \ra,
\Omega_{\La}^{\id \La}(M) )$. Therefore, as shown in the proof of
Theorem \ref{vangor1}, and by Theorem \ref{vangor2}, the implication
$$\Ext_{\La}^{i}( \Omega_{\La}^{\id \La}(M),C_N)=0 \text{
for } i \gg 0 \Rightarrow \Ext_{\La}^{i}( \Omega_{\La}^{\id
\La}(M),C_N)=0 \text{ for } i>0$$ holds, showing that (i) implies
(ii). To show that (iii) implies (iv), fix an exact sequence
$$0 \to Q_M \to C_M \to M \to 0$$
in $\mod \La$, in which $Q_M$ has finite projective (and injective)
dimension, and $C_M$ is maximal Cohen-Macaulay. There is an
isomorphism $\Ext_{\La}^i(N,M) \simeq \Ext_{\La}^i(N,C_M)$ for every
$i > \id \La$. Also, as above, since either $\Ext_{\La}^*(M, \La /
\ra )$ or $\Ext_{\La}^*( \La / \ra, M )$ belongs to $\Noethfl R$, so
does one of $\Ext_{\La}^*( C_M, \La / \ra )$ and $\Ext_{\La}^*( \La
/ \ra, C_M )$. Hence (iii) implies (iv) by Theorem \ref{vangor1} and
Theorem \ref{vangor2}.

By Theorem \ref{auscat}, the subcategory $\thick_{\underline{\MCM}(
\La )}(C_M)$ of $\underline{\MCM}( \La )$ is a left and right
Auslander subcategory. Moreover, by assumption $\underline{\MCM}(
\La )$ is weakly $d$-Calabi-Yau for some integer $d \in \mathbb{Z}$.
Therefore the implications
\begin{eqnarray*}
\Ext_{\La}^i(M,N)=0 \text{ for } i \gg 0 & \Leftrightarrow &
\Ext_{\La}^i(C_M,C_N)=0 \text{ for } i \gg 0 \\
& \Leftrightarrow & \Hom_{\underline{\MCM}( \La )}(C_M, \s^iC_N) =0
\text{ for } i \gg 0 \\
& \Leftrightarrow & \Hom_{\underline{\MCM}( \La )}(C_M, \s^iC_N) =0
\text{ for } i \in \mathbb{Z} \\
& \Leftrightarrow & \Hom_{\underline{\MCM}( \La )}(C_N, \s^iC_M) =0
\text{ for } i \in \mathbb{Z} \\
& \Leftrightarrow & \Hom_{\underline{\MCM}( \La )}(C_N, \s^iC_M) =0
\text{ for } i \gg 0 \\
& \Leftrightarrow & \Ext_{\La}^i(C_N,C_M)=0 \text{ for } i \gg 0 \\
& \Leftrightarrow & \Ext_{\La}^i(N,M)=0 \text{ for } i \gg 0
\end{eqnarray*}
hold, and the proof is complete.
\end{proof}

For an Artin algebra $\La$ with radical $\ra$, if $\Ext_{\La}^*( \La
/ \ra, \La / \ra ) \in \Noethfl R$ for some graded-commutative ring
$R$ acting centrally on $D^b ( \La )$, then $\Ext_{\La}^*(M,N) \in
\Noethfl R$ for all modules $M,N \in \mod \La$. Moreover, if this
holds, then $\La$ is automatically Gorenstein by \cite[Proposition
5.6]{BIKO}. Consequently, we obtain the following ``global version"
of Theorem \ref{localsym}.

\begin{theorem}\label{globalsym}
Let $\La$ be a stably symmetric Artin algebra with Jacobson radical
$\ra$, and suppose that $\Ext_{\La}^*( \La / \ra, \La / \ra )$
belongs to $\Noethfl R$ for some graded-commutative ring $R$ acting
centrally on $D^b ( \La )$. Then for all modules $M,N \in \mod \La$,
the following are equivalent:
\begin{enumerate}
\item[(i)] $\Ext_{\La}^i(M,N)=0$ for $i \gg 0$.
\item[(ii)] $\Ext_{\La}^i(M,N)=0$ for $i > \id \La$.
\item[(iii)] $\Ext_{\La}^i(N,M)=0$ for $i \gg 0$.
\item[(iv)] $\Ext_{\La}^i(N,M)=0$ for $i > \id \La$.
\end{enumerate}
\end{theorem}

Next, we include a special case of this theorem. Recall that for a
commutative Artin ring $k$, an Artin $k$-algebra $\La$ is
\emph{symmetric} if there is an isomorphism $\La \simeq \Hom_k( \La,
k)$ of $\La$-$\La$-bimodules. Such an algebra is necessarily
selfinjective.

\begin{corollary}\label{symmetricalg}
Let $\La$ be a symmetric Artin algebra with Jacobson radical $\ra$,
and suppose that $\Ext_{\La}^*( \La / \ra, \La / \ra )$ belongs to
$\Noethfl R$ for some graded-commutative ring $R$ acting centrally
on $D^b ( \La )$. Then for all modules $M,N \in \mod \La$, the
following are equivalent:
\begin{enumerate}
\item[(i)] $\Ext_{\La}^i(M,N)=0$ for $i \gg 0$.
\item[(ii)] $\Ext_{\La}^i(M,N)=0$ for $i >0$.
\item[(iii)] $\Ext_{\La}^i(N,M)=0$ for $i \gg 0$.
\item[(iv)] $\Ext_{\La}^i(N,M)=0$ for $i > 0$.
\end{enumerate}
\end{corollary}

\begin{proof}
By \cite[Corollary 4.4]{Mori}, a symmetric Artin algebra is stably
symmetric.
\end{proof}

We turn now to a particular class of algebras having finitely
generated cohomology in the sense of Theorem \ref{globalsym} and
Corollary \ref{symmetricalg}. Details concerning the following can
be found in \cite{SnashallSolberg} and \cite{Solberg}. Let $k$ be a
field and $\La$ a finite dimensional $k$-algebra, and denote the
enveloping algebra $\La \otimes_k \La^{\op}$ of $\La$ by $\Lae$. For
$n \ge 0$, the $n$th \emph{Hochschild cohomology} group of $\La$,
denoted $\HH^n ( \La )$, is the vector space $\Ext_{\Lae}^n( \La,
\La )$. The graded vector space $\HH^* ( \La ) = \Ext_{\Lae}^*( \La,
\La )$ is a graded-commutative ring with Yoneda product, and for
every $M \in \mod \La$ the tensor product $- \otimes_{\La} M$
induces a homomorphism
$$\HH^* ( \La ) \xrightarrow{\varphi_M} \Ext_{\La}^*(M,M)$$
of graded $k$-algebras. If $N \in \mod \La$ is another module and
$\eta \in \HH^* ( \La )$ and $\theta \in \Ext_{\La}^*(M,N)$ are
homogeneous elements, then the relation $\varphi_N( \eta ) \circ
\theta = (-1)^{|\eta||\theta|} \theta \circ \varphi_M( \eta )$
holds, where ``$\circ$" denotes the Yoneda product. Therefore the
Hochschild cohomology ring $\HH^* ( \La )$ acts centrally on $D^b(
\La )$.

Suppose now in addition that $\La$ is indecomposable as an algebra,
and that $\La / \ra \otimes_k \La / \ra$ is semisimple (as happens
for example when $k$ is algebraically closed). Furthermore, suppose
that $\La$ is a \emph{periodic algebra}. That is, there exists a
number $p>0$ such that $\La$ is isomorphic to $\Omega_{\Lae}^p ( \La
)$ as a left $\Lae$-module (i.e.\ as a bimodule). By
\cite{Erdmann2}, \cite{Erdmann3} and \cite{Erdmann4}, this happens
for example when $\La$ is a selfinjective Nakayama algebra, a
M\"obius algebra or a preprojective algebra (see also
\cite{ErdmannSkowronski}), and by \cite{Green} the condition implies
that $\La$ is selfinjective. Letting $Q_n$ denote the $n$th module
in the minimal projective $\Lae$-resolution of $\La$, we have an
exact sequence
$$0 \to \La \to Q_{p-1} \to \cdots \to Q_0 \to \La \to 0$$
of bimodules, and we denote this by $\mu$. This extension is an
element of $\HH^p ( \La )$. If $\theta$ is an element of $\HH^n (
\La )$ for some $n>p$, then $\theta = \bar{ \theta } \mu^i$ for some
$i$ and a homogeneous element $\bar{ \theta }$ of degree not more
than $p$. Hence the Hochschild cohomology ring $\HH^* ( \La )$ is
generated over $\HH^0 ( \La )$ by the finite set of $k$-generators
in $\HH^1 ( \La ), \dots, \HH^p ( \La )$, and therefore is
Noetherian. If $S$ is a simple non-projective $\La$-module, then
$\mu \otimes_{\La} S$ is the beginning of the minimal projective
resolution of $S$, since $\La / \ra \otimes_k \La / \ra$ is
semisimple. Therefore $S$ must be periodic with period dividing $p$.
If $N$ is any finitely generated $\La$-module and $\omega$ is an
element of $\Ext_{\La}^n ( S,N )$ for some $n>p$, then, as above,
$\omega = \bar{ \omega } ( \mu \otimes_{\La} S)$ for some element
$\bar{ \omega } \in \Ext_{\La}^m ( S,N )$ with $m \leq p$. Therefore
$\Ext_{\La}^* ( S,N )$ is finitely generated as a module over $\HH^*
( \La )$, and this shows that $\Ext_{\La}^*( \La / \ra, \La / \ra )$
is a finitely generated $\HH^* ( \La )$-module. The following result
is therefore an application of Corollary \ref{symmetricalg}.

\begin{theorem}\label{periodicalg}
Let $k$ be a field, let $\La$ be a symmetric periodic $k$-algebra
with Jacobson radical $\ra$, and suppose that $\La / \ra \otimes_k
\La / \ra$ is semisimple. Then for all modules $M,N \in \mod \La$,
the following are equivalent:
\begin{enumerate}
\item[(i)] $\Ext_{\La}^i(M,N)=0$ for $i \gg 0$.
\item[(ii)] $\Ext_{\La}^i(M,N)=0$ for $i >0$.
\item[(iii)] $\Ext_{\La}^i(N,M)=0$ for $i \gg 0$.
\item[(iv)] $\Ext_{\La}^i(N,M)=0$ for $i > 0$.
\end{enumerate}
\end{theorem}

We finish this paper with an example in which we look at
selfinjective Nakayama algebras. As we have seen, these algebras are
periodic and therefore have finitely generated cohomology. However,
the example shows that unless the algebra is symmetric, symmetry
does \emph{not} necessarily hold in the vanishing of cohomology.

\begin{example}
Let $\Gamma$ be the circular quiver
$$
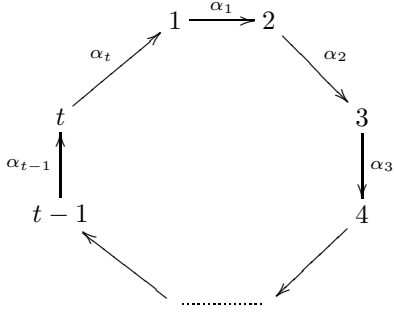
\xymatrix{
& 1 \ar[r]^{\alpha_1} & 2 \ar[dr]^{\alpha_2} \\
t \ar[ur]^{\alpha_t} &&& 3 \ar[d]^{\alpha_3} \\
t-1 \ar[u]^{\alpha_{t-1}} &&& 4 \ar[dl] \\
& \ar[ul] & \ar@{.}[l] }$$ where $t \geq 2$ is an integer. Let $k$
be a field, denote by $k \Gamma$ the path algebra of $\Gamma$ over
$k$, and let $J \subset k \Gamma$ be the ideal generated by the
arrows. Fix an integer $n \geq 1$, let $\La$ be the quotient algebra
$k \Gamma / J^{n+1}$, and denote by $\ra$ the Jacobson radical of
$\La$. Then $\La$ is a finite dimensional indecomposable
selfinjective Nakayama algebra, and $\Ext_{\La}^*( \La / \ra, \La /
\ra )$ is a finitely generated $\HH^* ( \La )$-module (the ring
structure of $\HH^* ( \La )$ was studied and determined in
\cite{Bardzell} and \cite{Erdmann2}).

Write $n=qt+r$, where $0 \leq r <t$. Let $S_i$ be the simple module
corresponding to the vertex $i$, and $P_i$ its projective cover.
There is an exact sequence
$$0 \to \Omega_{\La}^2 (S_i) \to P_{i+1 ( \mod t)}
\xrightarrow{\cdot \alpha_i} P_i \to S_i \to 0,$$ and it is easy
to see that $\Omega_{\La}^2 (S_i)$ is isomorphic to $S_{i+1+r(
\mod t )}$. Therefore the minimal projective resolution of $S_i$
is
$$\cdots \to P_{i+3+2r} \to P_{i+2+2r} \to P_{i+2+r} \to P_{i+1+r}
\to P_{i+1} \to P_i \to S_i \to 0,$$ with $\Omega_{\La}^{2j} (S_i) =
S_{i+j+jr}$ (all the indices are taken modulo $t$). A number of
completely different situations may occur, depending on the values
of the parameters $t$ and $r$. For example, if $r=0$, then we see
that all the simple modules appear infinitely many times as even
syzygies in the minimal projective resolution of any simple module.
Therefore, in this case, if $S$ and $S'$ are simple modules, then
$\Ext_{\La}^n(S,S')$ is nonzero for infinitely many $n$.

Note that when $r=0$, then $\La$ is symmetric, and so by Theorem
\ref{periodicalg} symmetry holds in the vanishing of $\Ext$.
However, symmetry does not hold for all Nakayama algebras. For
example, suppose $t \geq 3$ and $r=t-1$. Then the exact sequences
$$0 \to S_1 \to P_2 \to P_1 \to S_1 \to 0$$
$$0 \to S_2 \to P_3 \to P_2 \to S_2 \to 0$$
are the first parts of the minimal projective resolutions of $S_1$
and $S_2$, and therefore $\Ext_{\La}^n (S_1,S_2) \neq 0$ whenever
$n$ is odd, whereas $\Ext_{\La}^n (S_2,S_1) = 0$ for all $n$. Thus
in this situation there is no symmetry in the vanishing of $\Ext$
over $\La$.
\end{example}

\section*{Acknowledgements}

I would like to thank Steffen Oppermann and Idun Reiten for valuable
comments on this paper.

\end{document}